\documentclass[10pt,twocolumn,twoside]{IEEEtran}

\usepackage{graphicx}
\usepackage{array}
\usepackage{subfig}
\usepackage{cite}
\usepackage{stmaryrd}
\usepackage{amsfonts,latexsym,amssymb}
\usepackage[cmex10]{amsmath}
\usepackage{algorithm}
\usepackage{algorithmic}
\usepackage{hyperref}
\usepackage{multirow}
\usepackage{arydshln}

\input{zgx.math}

\newcommand{\nReal}{\Real_+}

\begin{document}
%
\title{Decomposition of Big Tensors With Low Multilinear Rank}
{\author{Guoxu~Zhou, Andrzej~Cichocki ~\IEEEmembership{Fellow,~IEEE}, and Shengli~Xie \IEEEmembership{Senior Member,~IEEE}

\thanks{Guoxu Zhou is with the Laboratory for Advanced Brain Signal Processing, RIKEN, Brain Science Institute, Wako-shi, Saitama 3510198, Japan. E-mail: zhouguoxu@ieee.org.}%
\thanks{Andrzej Cichocki is with the RIKEN BSI, Japan and Systems Research Institute, Warsaw Poland. E-mail: cia@brain.riken.jp.}%
\thanks{Shengli Xie is with the Faculty of Automation, Guangdong University of Technology, Guangzhou 510006, China. E-mail: eeoshlxie@scut.edu.cn.}}
}


\maketitle

\begin{abstract}
Tensor decompositions are promising tools for big data analytics as they bring multiple modes and aspects of data to a unified framework, which allows us to discover complex internal structures and correlations of data. Unfortunately most existing approaches are not designed to meet the major challenges posed by big data analytics. This paper attempts to improve the scalability of tensor decompositions and provides two  contributions: A flexible and fast algorithm for the CP decomposition (FFCP) of tensors based on their Tucker compression; A distributed randomized Tucker decomposition approach for arbitrarily big tensors but with relatively low multilinear rank. These two algorithms can deal with huge tensors, even if they are dense. Extensive simulations provide empirical evidence of the validity and efficiency of the proposed algorithms.
 \end{abstract}


%
\IEEEpeerreviewmaketitle

\section{Introduction}
\label{sec:Intro}


Tensor decompositions are emerging tools for analyzing multi-dimensional data that arise in many scientific and engineering fields such as signal processing, machine learning, chemometrics \cite{nmfbookCA,CPappSensor,FOBIUM2007,tensorComon}. Some data are naturally presented as tensors, such as color images, video clips, multichannel multi-trial EEG signals. In the era of big data, even more data with multiple aspects and/or modalities are collected and many of them are naturally correlated and can be organized as tensors for improved data analysis.  The methodologies that matricize the data and then apply matrix factorization approaches give \emph{flatten view} of data and often cause loss of internal structure information. Compared with matrix factorization techniques, tensor decomposition methodologies are versatile and provide multiple perspective stereoscopic view of data rather than a flatten one. Tensor decompositions are particularly promising for big data analysis as multiple data can be analyzed in a unified framework, which allows us to explore their  complex connections across multiple factor matrices simultaneously. By tensor decompositions, we achieve data compression, low-rank approximation, visualization, and feature extraction of high-dimensional multi-way data.

 \begin{table}[!t]
\caption{Notations and Definitions}
\label{tab:notations}
\centerline{
\begin{tabular}{ >{\hfill}p{.1\textwidth} | p{0.35\textwidth}  }
\hline \hline
\mat{A}, \mats[r]{a}, $a_{ir}$ & Matrix, the $r$th-column, and the ($i,r$)th-entry of matrix \mat{A}, respectively. \\
 $\Natural$, $\Real$  & Set of integers, real numbers.  $\Natural_+$ and $\nReal$ are their positive counterparts.  \\
$\Natural_N$  & The set of positive integers no larger than $N$ in ascending order, i.e., $\Natural_N=\set{1,2,\ldots,N}$. \\
$\tensor{Y}$, \tenmat{Y} & A  tensor, the mode-$n$ matricization of tensor \tensor{Y}. \\
$\hdp$, $\matdiv{}{}$  & Element-wise (Hadamard) product, division  of matrices or tensors. Moreover, we define $\frac{\tensor{A}}{\tensor{Y}}\defeq\matdiv{\tensor{A}}{\tensor{Y}}$.\\
$\kkp$, $\krp$ & Kronecker product and Khatri-Rao product (column-wise Kronecker product) of matrices. \\
\hline \hline
\end{tabular}
}
\end{table}

Two of the most popular tensor decomposition models are the canonical polyadic (CP) model, also known as APRAllel FACtor (PARAFAC) analysis, and the Tucker model, respectively \cite{Kolda09tensordecompositions}. The CP model, which is quite similar to the singular value decomposition (SVD) of matrices, decomposes the target tensor into the sum of rank-one tensors.  While CP decompositions (CPD) are able to give the most compact, essentially unique (under mild conditions \cite{Sidiropoulos2000}) representation of multiway data, they are unfortunately not always well-defined and may achieve a poor fit to the original data in practical applications. In contrast, Tucker decompositions estimate the subspace of each mode and the resulting representation of data are not unique, however, they often achieve more satisfactory fit. Hence these two models often serve different purposes.  Note that \apriori information on factor matrices could be incorporated in order to extracted physically meaningful components. For example, in \cite{SPL-SMBSS} statistical independence of components is used to improve the performance of CPD. Nonnegativity CPD, also called nonnegative tensor factorization (NTF) historically, is another important model that has attracted extensive study. In fact, nonnegative factorizations prove to be able to give useful parts-based sparse representation of objects and have been successfully applied in image interpretation and representation, classification, document clustering, etc \cite{nmfbookCA}. What is more, unlike unconstrained CPD, nonnegative CPD is always well defined and is free of ill-posedness of CPD \cite{NTFComon2009}.  

 \begin{figure}[!t]
\centerline{
    \includegraphics[width=\linewidth]{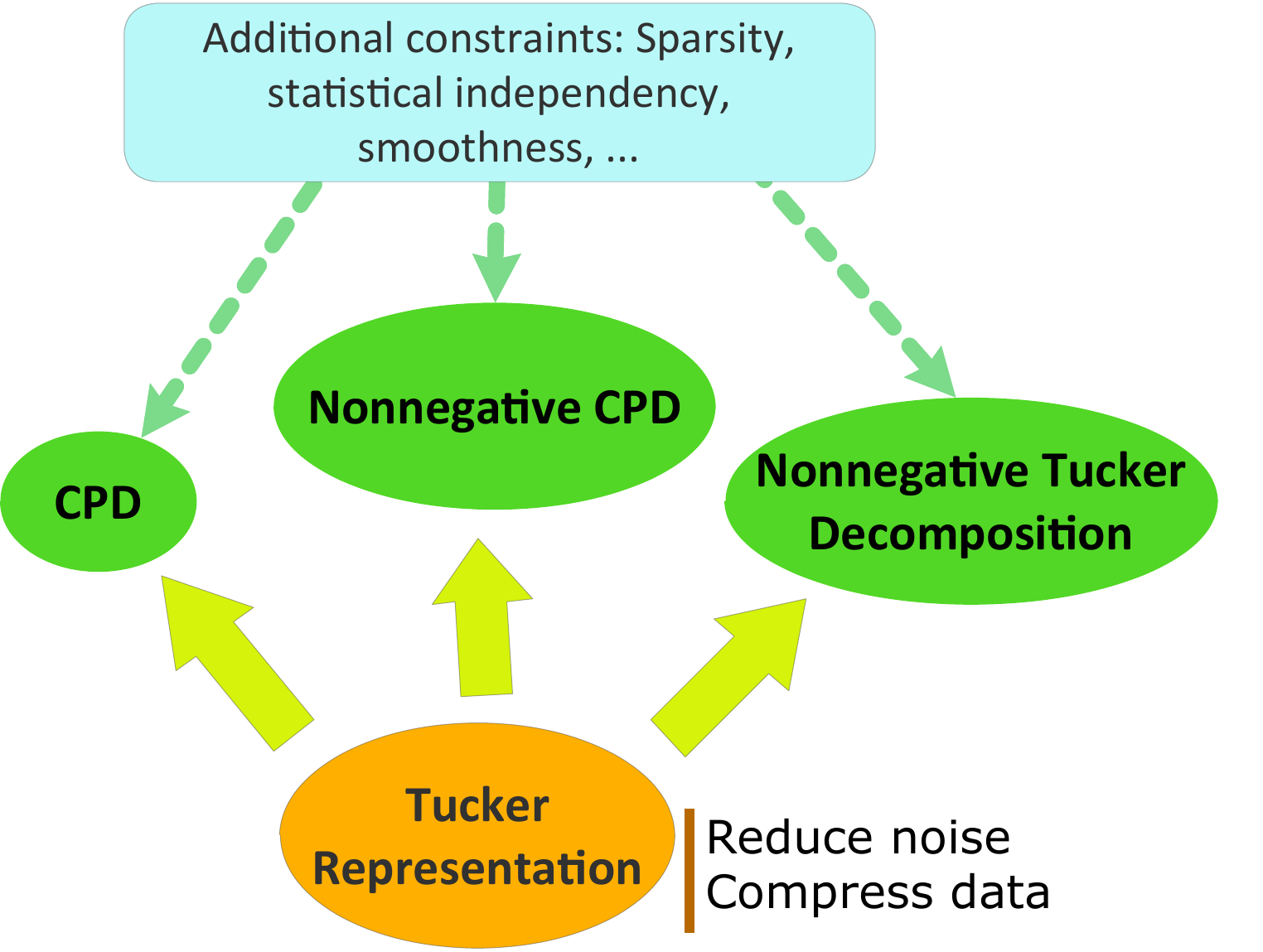}
}
\caption{Decompositions of a big tensor based on its compact Tucker representation.}
\label{fig:Tucker2AllTD}
\end{figure}

While many new algorithms for tensor decompositions have been developed in the last decade, only very recently has their scalability received special attention to meet the challenge posed by big tensor data analysis. For CP decompositions, the ParCube approach  \cite{ParCube} and the partition-and-merge approach in \cite{gridPARAFAC} were proposed, both of which depend on reliable CP decompositions of a set of small sub-tensors. However, existing CPD algorithms often suffer from the issue of ill-posedness in the sense that the optimal low-rank approximation does not exist at all. Hence in this type of methods, any ill-conditioned sub-tensors  formed by the sampling/partition operation  will unavoidably deteriorate  the final results. Another scalable approach GigaTensor \cite{GigaTensor} requires the target tensors are very sparse. For Tucker decompositions the MACH approach \cite{MACH} was also proposed to obtain a sparse tensor based on sampling the entries of a target tensor. Some other approaches that sample the columns (fibers) of target tensors can be found in \cite{randTenSVD,CURTensorCaiafaC2010}. However, these methods are often not suitable for very noisy data. Our purpose in this paper is to develop a general framework for big tensor decompositions. The basic idea is simple: Assuming the given big tensor has low multilinear rank, we use the Tucker model to compress the original big tensor first. Using this Tucker representation of data, both unconstrained/constrained CP and Tucker decompositions can be efficiently performed, as illustrated in \figurename \ref{fig:Tucker2AllTD}. We developed fast scalable algorithms for CPD based on a compact Tucker representation of data and  a new randomized algorithm for the Tucker decomposition of a big tensor. 

The rest of the paper is organized as follows. In Section \ref{sec:Models} the CP and Tucker decomposition models and algorithms are briefly overviewed. In Section \ref{sec:TPD} a family of new algorithms are proposed for (nonnegative) CP decompositions based on the Tucker representation of tensors. In Section \ref{sec:DRandTucker} two new distributed randomized Tucker decomposition algorithms are proposed that can be used to decompose big tensors with relatively low multilinear rank. Simulations results on synthetic and real-world data are presented in Section \ref{sec:simulations} and conclusions are made in Section \ref{sec:conclusion}.

In TABLE \ref{tab:notations} the frequently used notations are listed, and readers are referred to  \cite{Kolda09tensordecompositions, nmfbookCA} for more details on tensor conventions.

\section{Overview of Tucker and Polyadic Decompositions}
\label{sec:Models}
By using the Tucker model, a given high-order tensor $\tensor{Y}\in\Real^{I_1\times I_2\cdots\times I_N}$ is approximated as
\begin{equation}
\label{eq:Tucker}
\begin{split}
\tensor{Y}\approx\tensor{\hat{Y}}&=\tensor{G}\ttmn[1]{U}\ttmn[2]{U}\cdots\ttmn[N]{U}\\
&\defeq\compactTucker{G}{U}, 
\end{split}
\end{equation}
where $\matn{U}=\begin{bmatrix}
\matn{u}_1 & \matn{u}_2 & \cdots & \matn{u}_{R_n}
\end{bmatrix}\in\Real^{I_n\times R_n}$ is the mode-$n$ component (factor) matrix consisting of latent components $\matn{u}_{r_n}$ as its columns,  $\tensor{G}\in\Real^{R_1\times R_2\cdots\times R_N}$ is the core tensor reflecting the connections between the components, and  $(R_1,R_2,\ldots,R_N)$ is called the multilinear rank of \tensor{\hat{Y}} with $R_n=\rank{\tenmat{\hat{Y}}}$. 

By applying SVD to \tenmat{Y}, $n=1,2,\ldots,N$, we obtain the Tucker decomposition of a given tensor, which is referred to as the high-order SVD (HOSVD) and has been widely applied in multidimensional data analysis \cite{tensorFaces,HOgraphMatch,DenoisingHOSVD2013}. While the Tucker decomposition can be easily implemented, it is also criticized for the lack of uniqueness and suffering from curse of dimensionality. As another important tensor decomposition model, the CP model decomposes the target tensor into a set of rank-one terms, that is,
\begin{equation}
\label{eq:Polyadic}
\begin{split}
\tensor{Y}\approx\tensor{\hat{Y}} &=\sum_r \;\lambda_r\matn[1]{a}_r\outerp\matn[2]{a}_r\outerp\cdots\outerp\matn[N]{a}_r \\
&\defeq\compactcp{A},
\end{split}
\end{equation}
where $\outerp$ denotes the outer product of vectors, 
and, for simplicity, $\lambda_r$ is absorbed into $\matn[N]{a}_r$, $r=1,2,\ldots,R$.
CPD is essentially unique under mild conditions \cite{Sidiropoulos2000}, which has been found many important applications such as blind identification, feature extraction. However, the CPD of a given tensor is not always well-defined, in the sense that the optimal CP decomposition for a given rank may not exist at all \cite{deSilva:2008:ill-posted}.

\section{Efficient CP Decompositions of Tensors Based on Their Tucker Approximations}
\label{sec:TPD}

 \begin{figure*}[!t]
\centerline{
 \subfloat[Tucker+CP]{
    \includegraphics[width=0.48\linewidth]{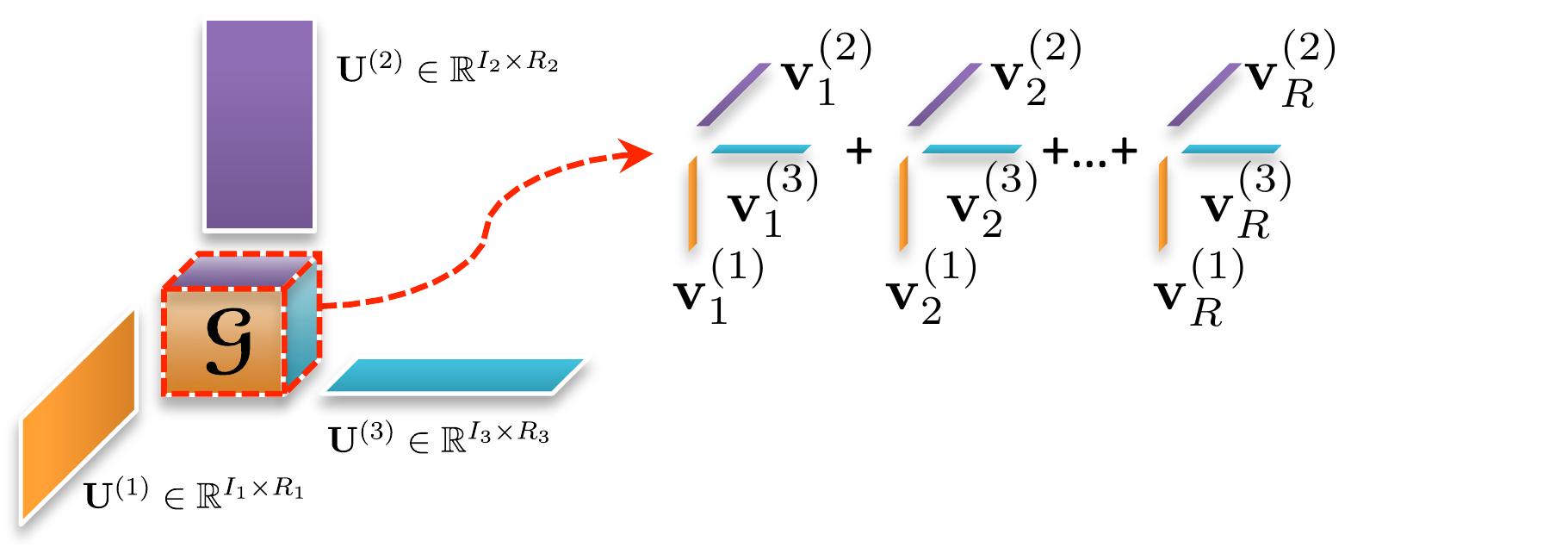}} \;
 \subfloat[FFCP]{
    \includegraphics[width=0.48\linewidth]{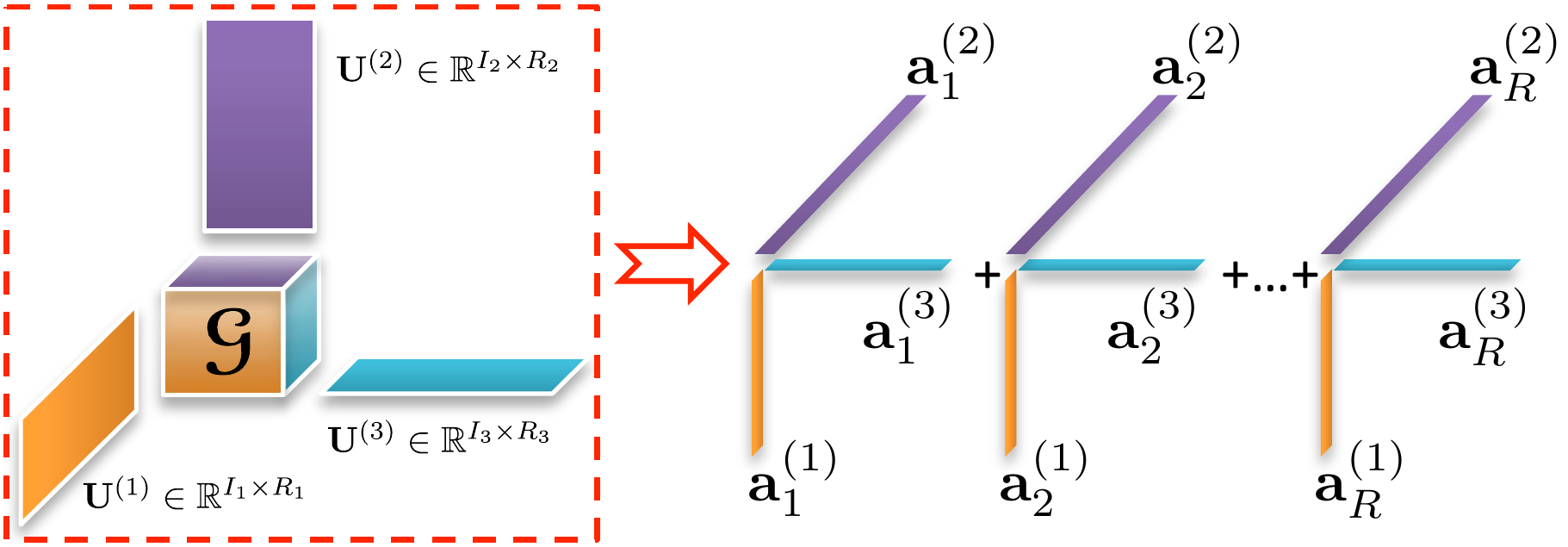}} 
}
\caption{Two different ways to perform CPD based on a Tucker approximation. In the Tucker+CP routine the CPD is performed on a small full tensor (i.e. the core tensor \tensor{G}), while in the FFCP method it is performed on a tensor in its approximate Tucker representation \compactTucker{G}{U}.}
\label{fig:TuckerCP}
\end{figure*}

\subsection{Overview of CPD Algorithms}
Alternating least-squares (ALS) has been a standard approach to perform CPD, which minimizes a total of $N$ LS problems 
\begin{equation}
\label{eq:PDLS}
J(\matn{A})=\frac{1}{2}\frob{\tenmat{Y}-\matn{A}\matn{B}{}^T}^2, n\in\Natural_N,
\end{equation}
in an alternating fashion in each iteration, possibly with additional constraints, e.g., nonnegativity, posed on the factor matrices \matn{A}, and 
 \begin{equation}
 \label{eq:PDBn}
 \matn{B}=\bigkrp\nolimits_{p\neq n}\matn[p]{A}\in\Real^{\prod_{p\neq n}I_p\times R}.  
 \end{equation}
Based on the gradient of $J(\matn{A})$
 \begin{equation}
 \label{eq:PDgrad}
 \fpartial{J}{\matn{A}}=\matn{A}\matn{B}{}^T\matn{B}-\tenmat{Y}\matn{B},
 \end{equation}
we can obtain the standard ALS update rule for the CPD
 \begin{equation}
 \label{eq:LS}
 \matn{A}\from\tenmat{Y}\matn{B}{(\matn{B}{}^T\matn{B})}^{-1},\; n=1,2,\ldots,N,
 \end{equation}
by letting $\fpartial{J}{\matn{A}}=0$; or a family of first-order nonnegative CPD algorithms \cite{SPM_NMFNTD}, for example, the one using the multiplicative update (MU) rule:
\begin{equation}
\label{eq:NPDMU}
\matn{A}\from\matn{A}\hdp\left[\matdiv{(\tenmat{Y}\matn{B})}{(\matn{A}\matn{B}{}^T\matn{B}})\right],
\end{equation} 
and the efficient hierarchical ALS (HALS) approach that updates only one column each time:
  \begin{equation}
\label{eq:HALS}
\matn{a}_r \from \matn{a}_r+\frac{1}{t^{(n)}_{r,r}}\proj_+{\left(\matn{q}_r-\matn{A}\matn{t}_r\right)},
\end{equation}
where $\proj_+(x)=\max(x,0)$ is element-wisely applied to a matrix/vector \cite{TSP-lraNMF,SPM_NMFNTD}, $t^{(n)}_{r,r}$ and $\matn{t}_r$ are the $(r,r)$ entry and the $r$th column of the matrix $\matn{T}\defeq\matn{B}{}^T\matn{B}$, respectively, and $\mats[r]{q}^{(n)}$ is the $r$th column of $\matn{Q}\defeq\tenmat{Y}\matn{B}$. Efficient computation of \eqref{eq:PDgrad} is the foundation of a large class of first-order methods for CPD and their variants \cite{SPM_NMFNTD}.   In \cite{spm2014_ConvOpti} the authors argued that first-order (gradient) methods generally provide acceptable convergence speed and are easy to parallelize, and hence more suitable than the second-order counterparts for large scale optimization problems. This motivates us to restrict our discussion to the first-order methods in this paper.

In the next we discuss how to compute the gradient terms in \eqref{eq:PDgrad} efficiently. While the term 
\begin{equation}
\label{eq:cpgradBtB}
\matn{B}{}^T\matn{B}=\bighdp\nolimits_{p\neq n}(\matn[p]{A}{}^T\matn[p]{A})\in\Real^{R\times R}
\end{equation}
in \eqref{eq:PDgrad} (and in \eqref{eq:LS}, \eqref{eq:NPDMU}, and \eqref{eq:HALS} as well)  can be efficiently computed, however, the term $\tenmat{Y}\matn{B}$ is quite computationally demanding, especially when \tensor{Y} is huge and $N$ is large. Indeed, it is likely that in many cases neither \tensor{Y} nor \matn{B} can be loaded into memory. To overcome this problem, we often hope to: 1) avoid accessing the huge tensor \tensor{Y} frequently (e.g., in every iteration); 2) avoid explicitly constructing the big matrix $\matn{B}\in\Real^{\prod_{k\neq n}I_k\times R}$. To this end, two properties of \tensor{Y} can be employed to improve scalability of CPD algorithms:
\begin{itemize}
\item
\emph{Sparsity of \tensor{Y}.} In this case the GigaTensor method \cite{GigaTensor} provides an elegant solution: only the nonzero entries of \tensor{Y} will be accessed and the columns of $\tenmat{Y}\matn{B}$ are computed in a parallel manner but without explicit construction of \matn{B};

\item
\emph{Low rank of \tensor{Y} (the tensor can be dense).  }  This type of methods use the low-rank approximation instead of the original huge dense tensor \tensor{Y}.  A standard method is to perform Tucker compression first,  typically by using the HOSVD, and then CPD is applied to the core tensor \cite{BroTuckerPD}, which is referred to as the Tucker+CP method in this paper. The approach in \cite{gridPARAFAC} splits the original big tensor into a set of small sub-tensors and then CPDs are applied to each sub-tensor. After that $\tenmat{Y}\matn{B}$ can be computed by using these sub-tensor in their CP approximation instead of \tensor{Y}. Arguing that some sub-tensors could be very ill-posed, which unavoidably causes poor performance, the PARACOMP method \cite{spm2014_paracomp} uses a set of random matrices that share a pre-specified number of common anchor rows to compress the tensor \tensor{Y} using the Tucker model and then a set of much smaller tensors are obtained. Then the CPD is performed on these small tensors simultaneously and original true components are recovered by solving a linear system incorporating the information provided by the anchor rows. This method avoids the SVD of huge matrices required by the HOSVD, however, it involves a total of $P$  Tucker compression procedures and CPD of $P$ small tensors, where  $P\ge {I_n}/{R}$.
\end{itemize}

Both the method in \cite{gridPARAFAC} and the PARACOMP method in \cite{spm2014_paracomp} rely on reliable CPD of all sub (or intermediate) tensors. However, it is widely known that existing CPD approaches often converge to local solutions, especially when bottlenecks exist in one or more factor matrices \cite{PComon2009} or even worse,  the CPD is simply not unique.  For example, in \cite{spm2014_paracomp} the authors argued that the method in  \cite{gridPARAFAC} assumed each sub tensor had unique CPD, which however cannot be guaranteed from the global uniqueness conditions alone \cite{spm2014_paracomp}. In the PARACOMP,  the uniqueness of CPD of each small tensor can be guaranteed from the global uniqueness conditions, however, it did not solve the local convergence problem and is sensitive to noise.  In other words, to achieve satisfactory results, the PARACOMP actually requires that all the CPDs in the intermediate step  converge globally, which is often unrealistic. Hence, we believe it is  important to develop more robust and scalable algorithms for the CPD, especially when \tensor{Y} is dense. In the meanwhile, we should avoid using CPDs as intermediate steps due to the local convergence issue.

Let us go into more detail of the Tucker+CP approach. In this approach, the big tensor \tensor{Y} is approximated by its HOSVD of $\tensor{Y}\approx\tensor{\tilde{Y}}=\compactTucker{G}{U}$, where $\matn{U}\in\Real^{I_n\times R_n}$ is computed from the truncated SVD of \tenmat{Y} with $\matn{U}{}^T\matn{U}=\matI$, then the CPD is applied to the much smaller tensor  $\tensor{G}\in\Real^{R_1\times R_2 \times\cdots\times R_N}$ such that
\begin{equation}
\label{eq:Tucker+CP}
\tensor{G}\approx\compactcp{V},
\end{equation}
thereby leading to the CPD of $\tensor{Y}\approx\compactcp{A}$ with $\matn{A}=\matn{U}\matn{V}$.
 For large-scale problems very often $I_n\gg R_n$ and the CPD of \eqref{eq:Tucker+CP} can be efficiently solved even using a standard desktop computer. 
 The Tucker+CP routine, which is originally proposed first  in \cite{BroTuckerPD}, is very efficient for the big tensors with low multilinear ranks. However, this framework is not directly applicable if  it is necessary  to impose additional constraints such as nonnegativity and sparsity on factor matrices. Very recently, Cohen, Farias and Comon proposed a novel efficient method to deal with the nonnegativity of factor matrices \cite{ComonNTF}. Our objective in this paper is to develop a more general and flexible scheme which can  not only achieve similar efficiency and scalability but also allow  to  impose  any  additional or optional constraints on factor matrices very easily.
 
\subsection{CPD of Big Tensors by Based on Their Tucker Compression}
The Tucker decomposition is suitable for pre-processing and compressing data because it is  significantly faster than the CPD and it has upper error bounds guarantee \cite{HOSVD2000}. In this subsection we consider a slightly modified version of the Tucker+CP approach to overcome its limitations. In the first step the Tucker compression is performed such that $\tensor{Y}\approx\tensor{\tilde{Y}}=\compactTucker{G}{U}$, as in the Tucker+CP approach. In the second step, however, we perform the CPD on the whole tensor \tensor{\tilde{Y}} instead of the core tensor $\tensor{G}$, by minimizing the following cost function:
\begin{equation}
\label{eq:FFCP}
\begin{split}
J(\matn{A}) & =\frac{1}{2}\frob{\tensor{\tilde{Y}}-\compactcp{A}}^2,\\
s.t. &\quad \matn{A}\in C_n,
\end{split}
\end{equation}
where $C_n$ denotes some constraints imposed on \matn{A}. 

\begin{algorithm}[!t]
\caption{The Flexible Fast Nonnegative CP (FFCP) Decomposition Algorithm}
\label{alg:FFCP}
\begin{algorithmic}[1]
 \REQUIRE $\tensor{Y}\in\Real^{I_1\times I_2 \cdots \times I_N}$, the output rank $R$.
 \STATE Perform Tucker decomposition to obtain $\tensor{Y}\approx\compactTucker{G}{U}$. Initialize \matn{A} and $\matn{V}=\matn{U}{}^T\matn{A}$.
 \REPEAT 
 \FOR {$n\in\Natural_N$}
 \STATE Efficiently compute $\matn{H}=\tenmat{G}\left[\bigkrp\nolimits_{p\neq n}{\matn[p]{V}}\right]$.
 \STATE Compute $\tenmat{Y}\matn{B}\approx\matn{U}\matn{H}$ and $\matn{B}{}^T\matn{B}=\bighdp_{p\neq n}(\matn[p]{A}{}^T\matn[p]{A})$. 
\STATE{ Update $\matn{A}$ by using \eqref{eq:LS}; or \eqref{eq:NPDMU}/\eqref{eq:HALS} to impose nonnegativity constraints; and/or \eqref{eq:Sparsity} to improve the sparsity of components.}
 \STATE {$\matn{V}\from\matn{U}{}^T\matn{A}$.}
 \ENDFOR
 \UNTIL{ a stopping criterion is satisfied}
\RETURN \tensor{G}\from\tensor{Y}, \matn{A}, $n\in\Set{N}$.
\end{algorithmic}
\end{algorithm}

With \eqref{eq:FFCP}, the term $\tenmat{Y}\matn{B}$ in the gradient $\fpartial{J}{\matn{A}}$ can be efficiently computed from: 
\begin{equation}
\label{eq:FFCPYB}
\tenmat{Y}\matn{B}\approx\matn{U}\tenmat{G}\left[\bigkrp\nolimits_{p\neq n}{\matn[p]{V}}\right],
\end{equation}
where
\begin{equation}
\label{eq:FFCPCmat}
\matn{V}=\matn{U}{}^T\matn{A}\in\Real^{R_n\times R}.
\end{equation}
Note that \matn{V} is quite small as $R_n\le R$, and  it can be computed by using distributed systems or the matrix partitioning method when $I_n$ is really huge. 
As such, the ALS update rule \eqref{eq:LS} can be replaced by
\begin{equation}
\label{eq:LRALS}
\matn{A}\from\matn{U}\tenmat{G}\left(\bigkrp_{p\neq n}{\matn[p]{V}}\right)\left(\bighdp_{p\neq n}\matn[p]{A}{}^T\matn[p]{A}\right)^{-1}.
\end{equation}
Because the approximation \tensor{\tilde{Y}} may contain negative entries, \eqref{eq:NPDMU} can be replaced by
\begin{equation}
\label{eq:LRANPD}
\matn{A}\from\matn{A}\hdp\frac{\proj_+\left(\matn{U}\tenmat{G}\left(\bigkrp\nolimits_{p\neq n}{\matn[p]{V}}\right)\right)}{\matn{A}\left(\bighdp_{p\neq n}\matn[p]{A}{}^T\matn[p]{A}\right)}.
\end{equation}
As a result, the nonnegativity of factor matrices can be easily incorporated in the proposed scheme. If we hope to extract sparse latent component, we consider the $l_1$-norm regularized CPD proposed in \cite{CP_L1} and minimize $J(\matn{A})+c\frob[1]{\matn{A}}$, thereby leading to
\begin{equation}
\label{eq:CPL1}
\matn{A}\from\mathcal{S}_{c}(\matn{A}),
\end{equation}
where the sparsity operator $\mathcal{S}_c(\cdot)$ is element-wisely defined as
\begin{equation}
\notag
\label{eq:Sparsity}
\mathcal{S}_c(x)=\begin{cases}
x-c; & \text{if}\; x>c; \\
0; & \text{if}\; -c\le x\le c; \\
x+c; & \text{if} \; x<-c.
\end{cases}
\end{equation}
 In summary, most first-order based constrained/unconstrained CPD algorithms can benefit from the proposed method. Approximately, the computational complexity is reduced from \bigO{R_n\prod_{p\in\Natural_N}I_p} in \eqref{eq:PDgrad} to only  about \bigO{I_nR_n^2+R_n\prod_{p\in\Natural_N}R_p}, due to using a compact representation of \tensor{Y}. 
 
 \begin{figure}[!t]
\centerline{
    \includegraphics[width=\linewidth]{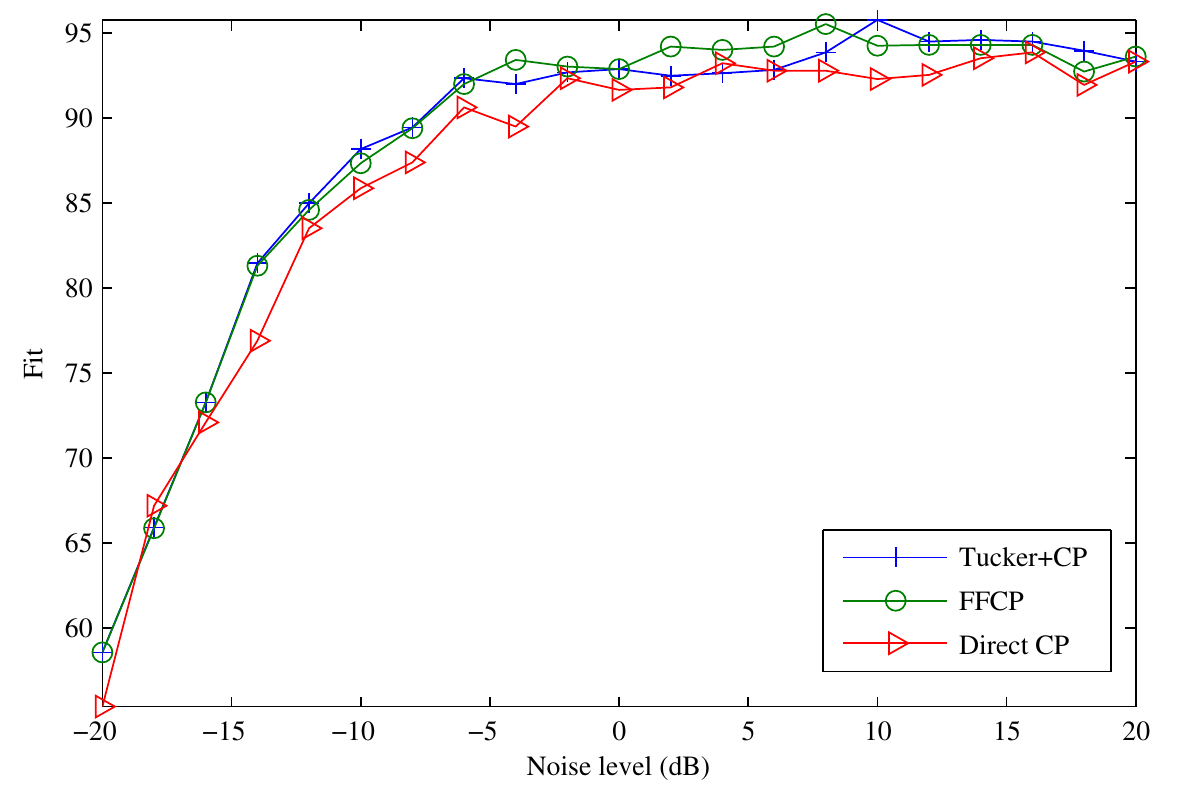}
}
\caption{The Fit values averaged over 50 Monte-Carlo runs obtained by 3 different ways (i.e., the Direct CP that applies CPD on \tensor{Y} directly, the FFCP on the Tucker format of $\tensor{Y}\approx\compactTucker{G}{U}$, and the Tucker+CP on the core tensor \tensor{G} first and then construct the CP approximation) to perform CPD on a synthetic tensor $\tensor{Y}\in\Real^{50\times50\times50\times50}$ with $R=10$. The data tensor \tensor{Y} was generated by using \eqref{eq:Polyadic} and different levels of Gaussian noise (SNR in dB) were added.}
\label{fig:3Ways}
\end{figure}

It is worth noticing that if $\tensor{\tilde{Y}}=\compactTucker{G}{U}$ with $\matn{U}\matn{U}{}^T=\matI$,  minimizing \eqref{eq:FFCP} is equivalent to solving:
\begin{equation}
\label{eq:FFCPvsCore}
\begin{split}
\min_{\matn{V}} \quad J(\matn{V})  & = \frac{1}{2} \frob{\tensor{G}-\compactcp{V}},\\
s.t. & \quad \matn{A}=\matn{U}\matn{V}\in C_n,
\end{split}
\end{equation}
where  \matn{V} is defined in \eqref{eq:FFCPCmat}, and $C_n$ denotes some constraints on \matn{A}. Indeed, \eqref{eq:FFCPvsCore} has been considered in \cite{ComonNTF}  specially for NTF. Although this method essentially has the same computational complexity as the FFCP method,  the main advantage of \eqref{eq:FFCP} is that any off-the-shelf (constrained) CPD methods can be directly applied with efficient computation of the gradient terms \eqref{eq:cpgradBtB} and \eqref{eq:FFCPYB},  without the need of designing a new algorithm to solve \eqref{eq:FFCPvsCore}.

In summary, not only does \eqref{eq:FFCP} benefit from the same advantage of low-computational complexity as the Tucker+CP method, it is also quite flexible because additional constraints on \matn{A} can be easily incorporated, such as nonnegativity and/or sparsity. For this reason, the new method is called the flexible fast CP (FFCP) decomposition algorithm and the corresponding pseudo code is briefly listed in Algorithm \ref{alg:FFCP}\footnote{ Although in Algorithm \ref{alg:FFCP} only the LS update and multiplicative update are listed, any update rules based on the gradient (e.g., HALS, see \cite{SPM_NMFNTD} for example) could be incorporated into the scheme of FFCP straightforwardly.}. In \figurename \ref{fig:3Ways} we compared the Tucker+CP, the FFCP (solving \eqref{eq:FFCP}), and the standard CPD (solving \eqref{eq:PDLS}) algorithms. The data tensor $\tensor{Y}\in\Real^{50\times50\times50\times50}$ was generated by using \eqref{eq:Polyadic} with $R=10$ and contaminated by different levels of noise (in SNR) ranging from -20dB to 20dB. In the figure we showed their Fit values averaged over 50 Monte-Carlo runs, where $\text{Fit}=1-{\frob{\tensor{Y}-\tensor{\hat{Y}}}}/{\frob{\tensor{Y}}}$. It turns out that the Tucker compression generally will not degenerate the accuracy. On the contrary, it may improve the stableness of CPD algorithms.


Although the matrices \matn{V} in \eqref{eq:FFCPYB} are generally small,  construction of the  tall $(\prod_{k\neq n}R_k)$-by-$R_n$ matrix, i.e., $\bigkrp_{k\neq n}\matn[p]{V}$, and  direct computation of $\matn{H}\defeq\tenmat{G}\left[\bigkrp\nolimits_{p\neq n}{\matn[p]{V}}\right]$ can be avoided by applying either one of the following methods:
\begin{itemize}
\item Notice that the $r$th column of \matn{H} can be computed using  $\matn{h}_r=\tensor{G}\times_{p\neq n}\matn[p]{a}_r{}^T$, which allows \matn{H} to be computed in parallel by using the shared data tensor \tensor{G} but without constructing  $\bigkrp_{k\neq n}\matn[p]{V}$ explicitly.

\item The term $\tenmat{G}\left[\bigkrp\nolimits_{p\neq n}{\matn[p]{V}}\right]$ can be computed in parallel by using the technique adopted in the GigaTensor \cite{GigaTensor}. Here \tenmat{G} is generally not sparse but relatively small.
\end{itemize}
Moreover, $\matn{H}\defeq\tenmat{G}\left[\bigkrp\nolimits_{p\neq n}{\matn[p]{V}}\right]$ can be viewed as the gradient term when we perform the CPD of $\tensor{G}\approx\compactcp{V}$ (see \eqref{eq:PDgrad} and \eqref{eq:FFCPvsCore}).  For this reason, most techniques which have been proposed to speed up standard CP decompositions  can be applied to compute \matn{H} straightforwardly. 
For very large-scale data where $I_n$ are extremely huge, the products of $\matn{A}{}^T\matn{A}$ and $\matn{U}{}^T\matn{A}$ can be performed on distributed systems, see, e.g., \cite{distMatProd2012}; or simply computed by using the matrix partitioning technique on a standard computer.

In summary, even with very limited computational resource, the FFCP algorithm can perform the (nonnegative) CPD of very large-scale tensors, provided that their low-rank Tucker representation can be obtained. Notice that the Tucker approximation of a given tensor plays a fundamental role in the proposed FFCP method as well as the Tucker+CP, the method proposed in \cite{ComonNTF}, and the PARACOMP \cite{paracomp}; and the HOSVD often serves as workhorse to perform Tucker compression. However, the HOSVD involves the eigenvalue decomposition of huge matrices \tenmat{Y}, which makes it unsuitable for large-scale problems \cite{paracomp}.  In the next section we focus on how to overcome this limitation.

\section{Distributed Randomized Tucker Decompositions (RandTucker)}
\label{sec:DRandTucker}

\subsection{Randomized Tucker Decomposition}
\label{sec:RandTucker}
In the Tucker decomposition of \eqref{eq:Tucker}, the factor matrices \matn{U} are estimated by minimizing the following LS problems:
\begin{equation}
\label{eq:TuckerLS}
J(\matn{U})=\frob{\tenmat{Y}-\matn{U}\matn{B}{}^T},
\end{equation}
where $\matn{B}=\left[\bigkkp\nolimits_{k\neq n}\matn[k]{U}\right]\tenmat{G}^T$.
In the HOSVD \cite{HOSVD2000} the factor matrix \matn{U} is estimated as the leading singular vectors of \tenmat{Y}. However, performing SVD is prohibitive for very large-scale problems.  Recently a set of probabilistic  algorithms have been developed to discover the low-rank structure of huge matrices \cite{siam_probLowRank}. Here we consider the following randomized range finder (RRF) algorithm (Algorithm 4.1 of \cite{siam_probLowRank}) to estimate the basis of \matn{U}:
\begin{enumerate}
\item Drawn a random Gaussian matrix $\mat{\Omega}\in\Real^{\prod_{p\neq n}I_p\times \tilde{R}_n}$.
\item Obtain the $I_n\times \tilde{R}_n$ matrix $\matn{Z}=\tenmat{Y}\matn{\Omega}$.
\item Compute $\matn{U}\in\Real^{I_n\times\tilde{R}_n}$ as an orthonormal basis of $\matn{Z}\in\Real^{I_n\times\tilde{R}_n}$ by using, e.g., QR decompositions.
\end{enumerate}
In the above $\tilde{R}_n=R_n+p$  and $p$ is referred to as the \emph{oversampling parameter}, and $p=5$ or $p=10$ was suggested in \cite{siam_probLowRank}. Based on this procedure, we propose a randomized Tucker decomposition (RandTucker) algorithm to perform Tucker decompositions, see Algorithm \ref{alg:RandTucker}. Compared with the HOSVD algorithm, the proposed RandTucker algorithm is highly parallalizable and suitable for distributed computation.

\begin{figure}[!t]
\centerline{
\includegraphics[width=\linewidth]{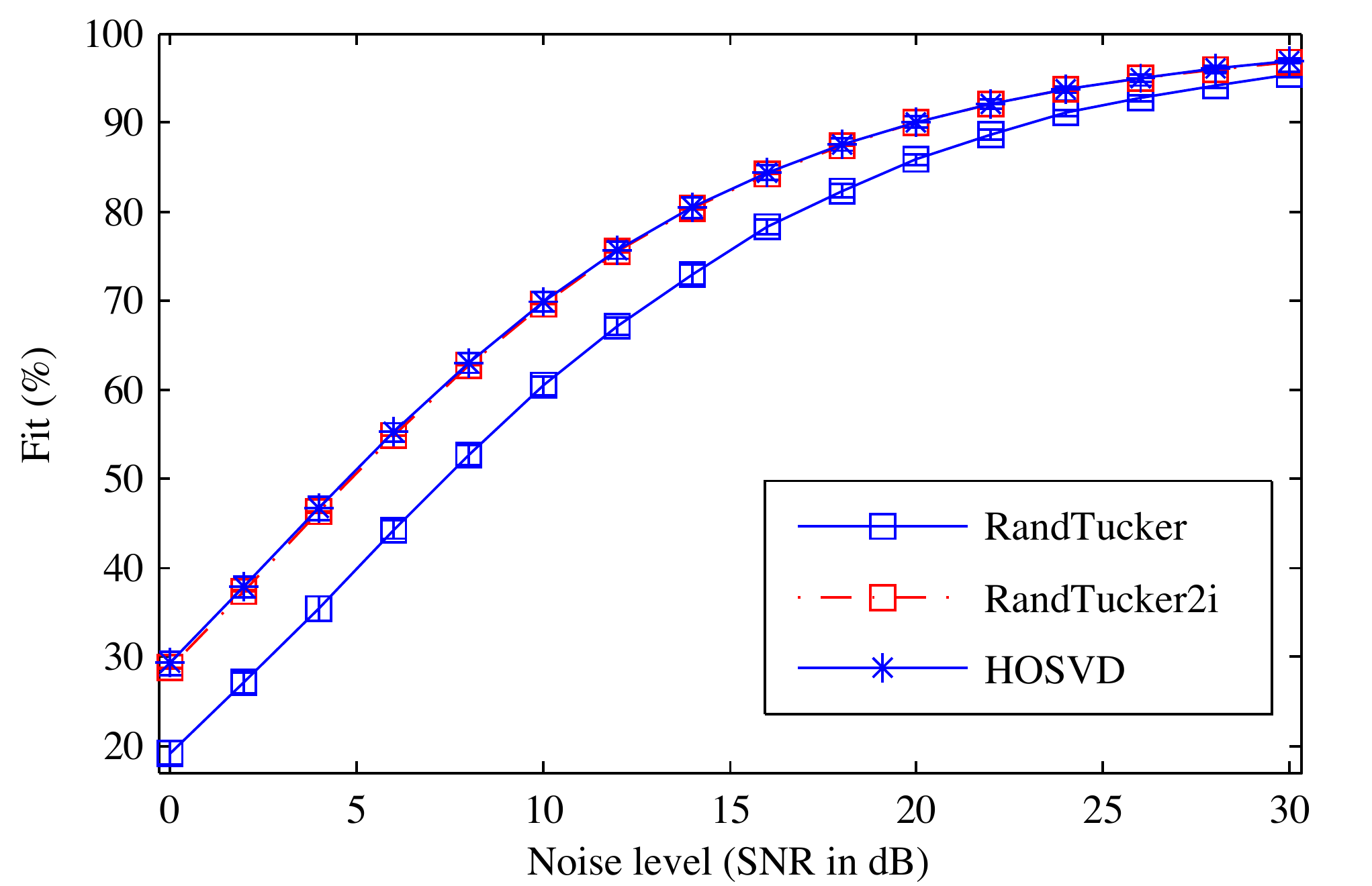}
}
\caption{Performance comparison between the RandTucker, the RandTucker2i, and the HOSVD algorithm in terms of mean value of Fit and the standard derivation over 50 Monte-Carlo runs.}
\label{fig:rtuckerhosvd}
\end{figure}

Several randomized algorithms have already been proposed for Tucker decompositions, such as the method using fiber sampling \cite{CURTensorCaiafaC2010,spm2014_Incomplete}. Compared with them, the proposed RandTucker method has an important advantage: it is able to obtain nearly optimal approximation even for very noisy data, by using only one additional ALS-like iteration, which is referred to as RandTucker2i, see Algorithm \ref{alg:RandTucker2i} and \figurename \ref{fig:rtuckerhosvd} for the performance comparison. In this figure, a 200-by-200-by-200 tensor \tensor{Y} was generated in each run by using the \eqref{eq:Tucker} with $R_n=10$, then different levels of Gaussian noise was added to \tensor{Y}  to generate the observation tensor. We plotted the Fit values (and the standard variations) obtained by the RandTucker (with $p=10$), the RandTucker2i ($p=10$), and the HOSVD averaged over 50 Monte-Carlo runs. From the figure, the RandTucker achieved slightly worse Fit values, especially when noise  was heavy, while the RandTucker2i achieved almost the same Fit values as the HOSVD. As the RandTucker and the RandTucker2i do not involve the singular value decomposition of big matrices and they are highly parallelizable, they are very promising in the Tucker decompositions of big tensors.

\begin{algorithm}[!t]
\caption{The Randomized Tucker Decomposition (RandTucker) Algorithm}
\label{alg:RandTucker}
\begin{algorithmic}[1]
 \REQUIRE $\tensor{Y}\in\Real^{I_1\times I_2 \cdots \times I_N}$ and the multilinear rank  $(\tilde{R}_1,\tilde{R}_2,\ldots,\tilde{R}_N)$.
\FOR{$n=1,2,\ldots,N$}
\STATE Compute  $\matn{Z}=\tenmat{Y}\matn{\Omega}$, where $\matn{\Omega}$ is an $(\prod_{k\neq n}I_k)$-by-$\tilde{R}_n$ random Gaussian matrix.
\STATE Compute \matn{U} as an orthonormal basis of \matn{Z} by using, e.g, the QR decomposition.
\STATE $\tensor{Y}\from\tensor{Y}\ttmn{U}{}^T$, $I_n\from \tilde{R}_n$.
\ENDFOR
\STATE{\tensor{G}\from\tensor{Y}.}
\RETURN $\tensor{Y}\approx\compactTucker{G}{U}$.
\end{algorithmic}
\end{algorithm}

\begin{algorithm}[!t]
\caption{The RandTucker Algorithm with 2 Iterations only (RandTucker2i)}
\label{alg:RandTucker2i}
\begin{algorithmic}[1]
 \REQUIRE $\tensor{Y}\in\Real^{I_1\times I_2 \cdots \times I_N}$ and the multilinear rank  $(\tilde{R}_1,\tilde{R}_2,\ldots,\tilde{R}_N)$.
 \STATE{Initialize $\matn{U}\in\Real^{I_n\times \tilde{R}_n}$ as random Gaussian matrices.}
\STATE{\%\emph{ Repeat lines (3)-(7) twice.}}
\FOR{$n=1,2,\ldots,N$}
\STATE{$\tensor{X}=\tensor{Y}\times_{p\neq n}\matn[p]{U}$.}
\STATE Compute  $\matn{Z}=\tenmat{X}\matn{\Omega}$, where $\matn{\Omega}$ is an $(\prod_{k\neq n}\tilde{R}_k)$-by-$\tilde{R}_n$ random Gaussian matrix.
\STATE Compute \matn{U} as an orthonormal basis of \matn{Z} by using, e.g., the QR decomposition.
\ENDFOR
\STATE{$\tensor{G}\from\tensor{Y}\ttmn[1]{U}\ttmn[2]{U}\cdots\ttmn[N]{U}$}.
\RETURN $\tensor{Y}\approx\compactTucker{G}{U}$.
\end{algorithmic}
\end{algorithm}


 \begin{figure*}[!t]
    \centerline{
    \subfloat[]{
    \includegraphics[width=0.45\textwidth]{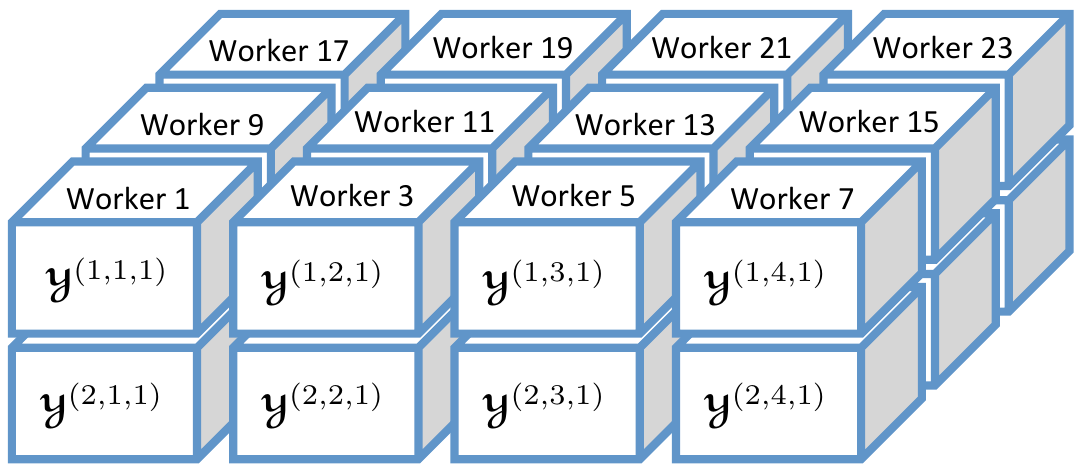}} \hfill
        \subfloat[]{
    \includegraphics[width=0.45\textwidth]{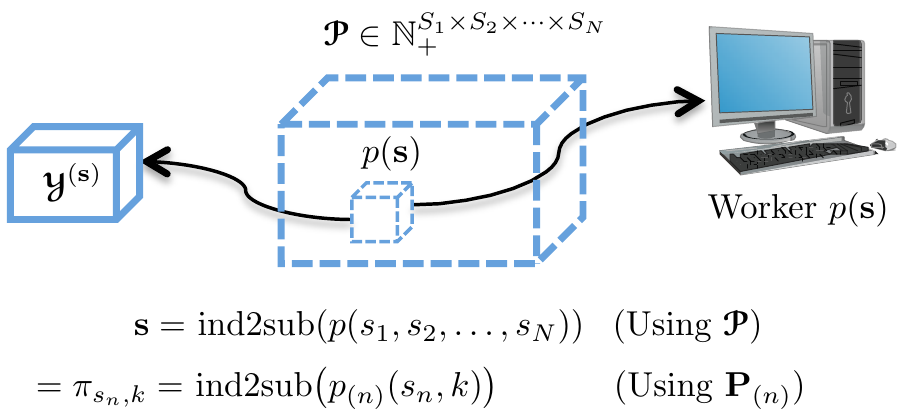}}
    }
    \caption{ Demonstration of small data tensors and the virtual computing resource tensor \tensor{P} in a cluster. 
    (a) An example: $\tensor{Y}\in\Real^{2I\times 4J \times 3K}$ is partitioned into smaller blocks $\tensor{Y}^{(\mat{s})}$ with the same size $I\times J \times K$, thereby leading to the $2\times4\times3$ computing resource tensor $\tensor{P}$ identifying a total of 24 workers. 
    (b) An entry $p(\mat{s})\in\set{1,2,\ldots,\prod_nS_n}$ of the computing resource tensor \tensor{P} indicates that Worker $p(\mat{s})$ processes the sub-tensor $\tensor{Y}^{(\mat{s})}$, where $\mat{s}=[s_1,s_2,\ldots,s_N]^T=\text{ind2sub}(p(\mat{s}))$.   Hence,  in (a), $p(1,3,2)=13$ means Worker 13 processes the sub-tensor $\tensor{Y}^{(1,3,2)}$, as $[1,3,2]=\text{ind2sub}(13)=p_{(2)}(3,3)$.
    }
    \label{fig:PartionOfTensor}
  \end{figure*}
\subsection{Distributed Randomized Tucker Decomposition}
\label{sec:DisTucker}
In the above RandTucker and RandTucker2i algorithms we have assumed that all the data can be processed in physical memory directly. For big data this could be unrealistic. A common scenario is that only partial data can be loaded in memory. In this subsection we discuss how to solve this problem  by using a computer cluster (in this paper we used the MATLAB Distributed Computing Server\footnote{See \url{www.mathworks.com/help/pdf_doc/mdce/mdce.pdf}.}). Here we consider a cluster which is composed of
\begin{itemize}
\item a number of workers with identity labels, each of which is a computational engine that executes computational tasks simultaneously. In the meanwhile they may communicate with each other and  the client as well;

\item a client which is used to schedule tasks to workers and gather results from workers, and return final results.
\end{itemize}

Given a big $N$th-order tensor $\tensor{Y}\in\Real^{I_1\times I_2\times \cdots\times I_N}$, we partition each dimension $I_n$ into a total of $S_n$ disjoint subsets, each of which consists of $L_{S_n}$ indices with $\sum_{s_n=1}^{S_n}L_{s_n}=I_n$. As such, the original big tensor is partitioned into a total of $W=S_1S_2\cdots S_N$ sub-tensors, and each sub-tensor can be denoted by $\tensor{Y}^{(s_1,s_2,\ldots,s_N)}$ with the size of ${L_{s_1}\times L_{s_2}\times\cdots\times L_{s_N}}$, where $s_n\in\Natural_{S_n}$, $n=1,2,\ldots,N$. At the same time suppose we have a computer cluster with $W$ workers (computational engines)\footnote{We have assumed that each sub-tensor can be efficiently handled in by the associated worker, otherwise we increase the number of workers and decrease the size of each sub-tensor accordingly}. We imagine that the $W$ workers are virtually spatially  permuted to form a $S_1\times S_2\times\cdots\times S_N$  tensor such that each worker is just associated with one sub-tensor, see \figurename \ref{fig:PartionOfTensor}(a). 
In other words, we will have two types of tensors:
\begin{itemize}
\item The small data tensors obtained by partitioning \tensor{Y}: 
\begin{equation}
\label{eq:subtensorBs}
\tensor{Y}^{(\mat{s})}=\tensor{Y}^{(s_1,s_2,\ldots,s_N)}\in\Real^{L_{s_1}\times L_{s_2} \times \cdots \times L_{s_N}}
\end{equation}
where $\mat{s}=\begin{bmatrix}
s_1 & s_2 & \cdots & s_N
\end{bmatrix}^T$ is used to specify the $(s_1,s_2,\ldots,s_N)$th-block of \tensor{Y}, $s_n\in\Natural_{S_n}$, $n=1,2,\ldots,N$; see \figurename \ref{fig:PartionOfTensor}(a) for an illustration.

\item
The computing resource tensor $\tensor{P}\in\Natural_+^{S_1\times S_2 \times \cdots \times S_N}$ formed by reshaping $[1,2,\ldots,W]$ that identifies the $W=S_1S_2\cdots S_N$ workers and the associated data tensors. The entries of \tensor{P} have a one-to-one correspondence with the sub-tensors $\tensor{Y}^{(\mat{s})}$, meaning that \emph{ Worker $p(\mat{s})$ process the sub-tensor $\tensor{Y}^{(\mat{s})}$, see \figurename \ref{fig:PartionOfTensor}(b).}


\end{itemize}
 It can be seen the RandTucker and the RandTucker2i majorly use four basic operations: 1)  Mode-$n$ unfolding of a tensor; 2) Matrix-matrix product to obtain $\tenmat{Y}\matn{\Omega}$;  3) Generating a Gaussian random matrix \mat{\Omega}; 4) Estimating an orthogonal basis of the range of $\tenmat{Y}\matn{\Omega}$. Below we discuss how to perform these four operations in a distributed computing system.

{\bf Mode-$n$ unfolding of tensor \tensor{Y}.} 
We unfold the computing resource tensor $\tensor{P}\in\Natural^{S_1\times S_2\times\cdots\times S_N}_+$ along mode-$n$ to obtain $\tenmat{P}\in\Natural^{S_n\times K}_+$, $K=\prod_{j\neq n}S_j$. Let $p_{(n)}(s_n,k)$  be the $(s_n,k)$th entry of $\tenmat{P}$, and define
\begin{equation}
\label{eq:ind2sub}
\begin{split}
& \mat{s} =\boldsymbol{\pi}_{s_n,k}=\text{ind2sub}(p_{(n)}(s_n,k)),\\
& s_n=1,2,\ldots,S_n, \quad k=1,2,\ldots,\prod\nolimits_{j\neq n}S_j,
\end{split}
\end{equation}
where the function $\text{ind2sub}(\cdot)$ converts  an entry of \tenmat{P}, i.e. $p_{(n)}(s_n,k)$, a worker identifier, into an $N$-by-1 vector of subscript values determining the position  of $p_{(n)}(s_n,k)$ in the tensor \tensor{P}, which allows us to locate the sub-tensor associated with this worker. In other words, the above \eqref{eq:ind2sub} means that the tensor $\tensor{Y}^{\pi_{s_n,k}}$  is processed by Worker $p_{(n)}(s_n,k)$. Then it can be seen that
\begin{equation}
\label{eq:YnGrid}
\tenmat{{Y}}=\begin{bmatrix}
\tenmat{Y}^{\boldsymbol{\pi}_{1,1}} & \tenmat{Y}^{\boldsymbol{\pi}_{1,2}} & \ldots & \tenmat{Y}^{\boldsymbol{\pi}_{1,K} } \\
\tenmat{Y}^{\boldsymbol{\pi}_{2,1}} & \tenmat{Y}^{\boldsymbol{\pi}_{2,2}} & \ldots & \tenmat{Y}^{\boldsymbol{\pi}_{2,K} } \\
\ldots & \ldots & \ldots & \ldots \\
\tenmat{Y}^{\boldsymbol{\pi}_{S_n,1}} & \tenmat{Y}^{\boldsymbol{\pi}_{S_n,2}} & \ldots & \tenmat{Y}^{\boldsymbol{\pi}_{S_n,K} } \\
\end{bmatrix}\in\Real^{I_n\times\prod_{j\neq n}I_j},
\end{equation}
where $\tenmat{Y}^{\boldsymbol{\pi}_{s_n,k}}$ is the mode-$n$ martricization of $\tensor{Y}^{\boldsymbol{\pi}_{s_n,k}}$ processed by  Worker $p_{(n)}(s_n,k)$, $s_n\in\Natural_{S_n}$, and $k\in\Natural_K$. 

\emph{In summary, to obtain \tenmat{Y}, we  perform mode-$n$ unfolding of $\tenmat{Y}^{(\mat{s})}$ simultaneously. 
Then the ($s_n,k$)th sub-matrix of \tenmat{Y} will be given by Worker $p_{(n)}(s_n,k)$, as shown in \eqref{eq:YnGrid}.}

 \begin{figure*}[!t]
    \centerline{
    \subfloat[]{
    \includegraphics[width=0.48\textwidth,height=0.32\textwidth]{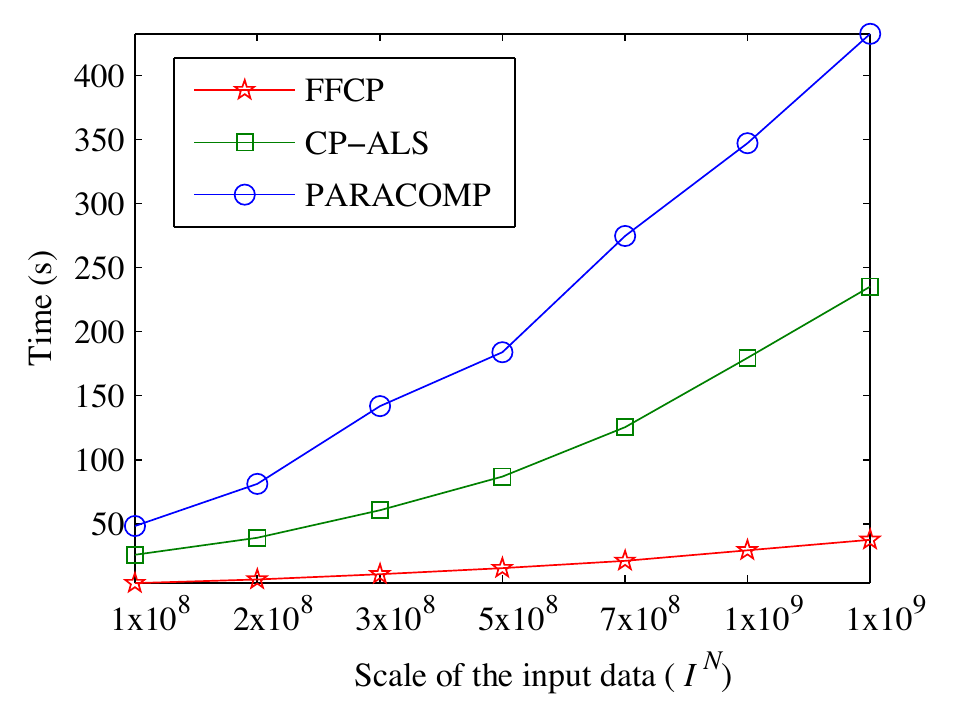}}\hfil
    \subfloat[]{
    \includegraphics[width=0.48\textwidth,height=0.32\textwidth]{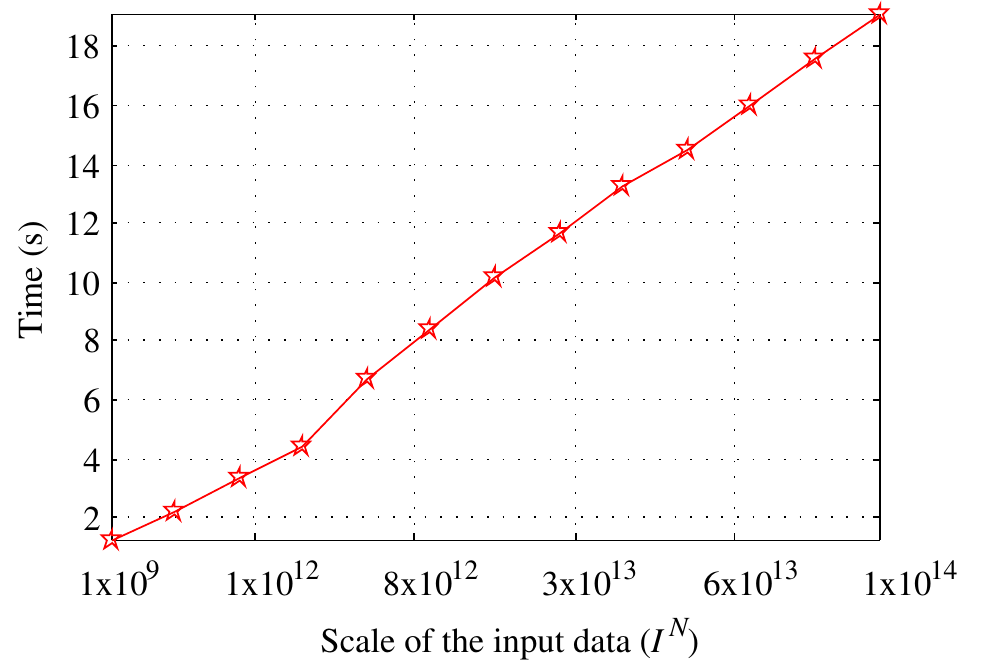}}
    }
    \caption{(a) Time consumption comparison of the FFCP, the CP-ALS, and the PARACOMP on condition that the data tensor of $I\times I\times I$ can be loaded into memory. For the FFCP the RandTucker2i was used to compress the data with $p=15$. (b) Scalability of the FFCP algorithm provided that the given big tensor was in its Tucker format. The consumed time by the FFCP increases nearly linearly with respect to $I$.}
    \label{fig:scalability}
  \end{figure*}

{\bf Matrix-matrix product.} Large scale matrix-matrix product can be performed by using any existing technique. As here the mode-$n$ unfolding matrix \tenmat{Y} has already been partitioned and distributed, we apply the matrix partitioning technique: we compute the product $\tenmat{Y}\matn{\Omega}$ by using
\begin{equation}
\label{eq:distYnOmega}
\begin{bmatrix}
\matn{Z}_1 \\
\matn{Z}_2 \\
\cdots \\
\matn{Z}_{S_n}
\end{bmatrix}=
\tenmat{Y}\matn{\Omega}=\tenmat{Y}\begin{bmatrix}
\matn{\Omega}_1 \\ \matn{\Omega}_2 \\ \ldots \\\matn{\Omega}_K 
\end{bmatrix}
\in\Real^{I_n\times \tilde{R}_n},
\end{equation}
where \tenmat{Y} is given by \eqref{eq:YnGrid}, $\matn{\Omega}_k\in\Real^{\prod_{j\neq n}L_{S_j}\times \tilde{R}_n}$ are the sub-matrices by partitioning the random Gaussian matrix \matn{\Omega} into $K$ parts. From \eqref{eq:distYnOmega}, we have $\matn{Z}_{s_n}=\sum_{k=1}^K\tenmat{Y}^{\boldsymbol{\pi}_{sn,k}}\matn{\Omega}_k$, and all the products $\tenmat{Y}^{\boldsymbol{\pi}_{s_n,k}}\matn{\Omega}_k$ are computed on Workers $p(\boldsymbol{\pi}_{s_n,k})$ simultaneously and then  are broadcasted to  Worker $p_{(n)}(s_n,k_0)$, typically $k_0=1$. Finally, for each $s_n\in\Natural_{S_n}$, we sum $\tenmat{Y}^{\boldsymbol{\pi}_{s_n,k}}\matn{\Omega}_k$ up over $k$ on Worker  $p_{(n)}(s_n,k_0)$ to yield $\matn{Z}_{s_n}$.

{\bf Generating the Gaussian random matrix \matn{\Omega}.} Note that the workers identified by the $k$th column of \tenmat{P}, i.e., the  workers $p_{(n)}(s_n,k)$, $s_n=1,2,\ldots, S_n$, share a same random Gaussian matrix $\matn{\Omega}_k$ (corresponding to the $k$th column of blocks of \tenmat{Y} in \eqref{eq:YnGrid}).  Hence $\matn{\Omega}_k$ can be generated on one worker, say $p_{(n)}(1,k)$, and then is broadcasted to all the workers $p_{(n)}(s_n,k)$ with $s_n\neq 1$. However, this way  suffers from heavy communication overhead. A more economical method is to broadcast the parameters (or the random number generators) that control the random number generation and then the random matrices $\matn{\Omega}_k$ are generated on all workers by using the received random number generators, as such, for each $k$, the workers $p_{(n)}(s_n,k)$, $s_n\in\Natural_{S_n}$, will generate identical random Gaussian matrices $\matn{\Omega}_k$ simultaneously. 

{\bf Estimating an orthonormal basis of the columns of $\tenmat{Y}\matn{\Omega}$}. To this end, we perform eigenvalue decomposition of $\sum_{s_n=1}^{S_n}\matn{Z}_{s_n}{}^T\matn{Z}_{s_n}=\mat{U\Sigma}\mat{U}^T$. Then the mode-$n$ factor matrix \matn{U} with orthonormal columns can be computed from
\begin{equation}
\label{eq:DRTAn}
 \matn{U}_{s_n}=\matn{Z}_{s_n}\mat{U}\mat{\Sigma}^{-\frac{1}{2}}
\end{equation}
in parallel on the workers $p_{(n)}(s_n,1)$, where $\matn{U}_{s_n}$ is the $s_n$th sub-matrix of \matn{U}, $s_n\in\Natural_{S_n}$.

After \matn{U} has been estimated, \tensor{Y} will be updated as $\tensor{Y}\from\tensor{Y}\ttmn{{U}}{}^T$, which can be realized  using the above mode-$n$ tensor unfolding and matrix-matrix product. After update the $n$th dimension of $\tensor{Y}^{(\mat{s})}$ will be reduced to be $\tilde{R}_n$ from $L_{S_n}$.

We repeat the above procedure for $n=1,2,\ldots,N$. As the last step, the core tensor is computed as $\tensor{G}=\sum_{w=1}^W\tensor{Y}^{(w)}$, i.e. the sum of sub-tensors distributed on all workers.

\subsection{Nonnegative Tucker Decomposition}
In our recent papers \cite{SPM_NMFNTD, lraNTD} we also showed how to efficiently perform nonnegative Tucker decomposition based on a Tucker approximation. The idea of that method is similar to the FFCP method: after the original high-dimensional raw data has been replaced by its Tucker approximation with low multilinear rank, the computational complexity of gradients with respect to the nonnegative core tensor and factor matrices can be significantly reduced in terms of both time and space. Readers are referred to \cite{SPM_NMFNTD, lraNTD}  for more details. In summary, Tucker compression can play a fundamental role to achieve highly efficient tensor decompositions, as illustrated in \figurename \ref{fig:Tucker2AllTD}. Based on the proposed randomized Tucker decomposition methods, we can easily perform various types of constrained/unconstrained tensor decompositions for very large-scale tensor data.



\begin{table*}
  \caption{Comparison between the algorithms when they were applied to perform CPD of synthetic data. The noises were Gaussian with SNR=10dB in all the cases.}
  \label{tab:simuCP}
  \centering
  \begin{tabular}{ l || c c  |  c c | c c}
  \hline \hline
  \multirow{2}{*}{Algorithm}
  & \multicolumn{2}{c|}{$N=3,\;I=200,\; R=10$}   & \multicolumn{2}{c |}{$N=3,\;I=500,\; R=20$} & \multicolumn{2}{c}{$N=6,\; I=20,\; R=5$}   \\
  \cline{2-7}
  & Fit (\%) & Time (s)    & Fit (\%) & Time (s)    & Fit (\%) & Time (s) \\ \cline{1-7} 
CP-ALS & $82.1\pm{10.2}$& 39.5 & $80.5\pm{6.8}$& 64.0 & $89.6\pm{11.2}$& 18.8     \\
FastCPALS & $79.8\pm{9.0}$& 10.4 & ${\bf82.9\pm{5.6}}$& 21.1 & $92.4\pm{10.9}$& 6.6     \\
PARAFAC & ${\bf82.5\pm{7.9}}$& 78.2 & $80.5\pm{6.1}$& 158.3 & $87.6\pm{11.4}$& 64.3     \\
FFCP & $81.9\pm{8.0}$& {\bf3.5} & $82.8\pm{5.6}$& {\bf4.0} & ${\bf93.4\pm{8.6}}$& {\bf5.0}     \\
\hline \hline
  \end{tabular}
\end{table*}

\begin{table*}
  \caption{Comparison between the algorithms when they were applied to perform nonnegative CPD of synthetic data. The noises were Gaussian with SNR=10dB in all the cases.}
  \label{tab:simuNTF}
  \centering
  \begin{tabular}{ l || c c  |  c c | c c}
  \hline \hline
  \multirow{2}{*}{Algorithm}
  & \multicolumn{2}{c|}{$N=3,\;I=200,\; R=10$}   & \multicolumn{2}{c |}{$N=3,\;I=500,\; R=20$} & \multicolumn{2}{c}{$N=6,\; I=20,\; R=5$}   \\
  \cline{2-7}
  & Fit (\%) & Time (s)    & Fit (\%) & Time (s)    & Fit (\%) & Time (s) \\ \cline{1-7} 

CP-NMU & $97.3\pm{1.2}$& 66.4 & $97.6\pm{0.8}$& 1695.2 & $91.7\pm{5.6}$& 426.5     \\
CP-HALS & $97.8\pm{1.0}$& 12.1 & $97.7\pm{0.7}$& 228.1 & ${\bf95.1\pm{4.1}}$& 101.4     \\
FFCP-MU & $99.0\pm{0.5}$& {\bf2.3} & $99.3\pm{0.2}$& {\bf12.1} & $93.4\pm{4.1}$& {\bf7.0}     \\
FFCP-HALS & ${\bf99.1\pm{0.1}}$& 6.1 & ${\bf99.5\pm{0.1}}$& 29.9 & $94.3\pm{3.5}$& 11.8     \\
  \hline \hline
  \end{tabular}
\end{table*}

\section{Simulations}
\label{sec:simulations}
\subsection{Synthetic data}
In this subsection we investigate the performance of the proposed FFCP algorithm by using synthetic data. All the experiments were performed in 64-bit MATLAB on a computer with Intel Core i7 CPU (3.33GHz) and 24GB memory, running 64-bit Windows 7. 

{\bf Simulations on CPD.} In each run a tensor $\tensor{Y}^*$ was generated by using the model \eqref{eq:Polyadic}, where the entries of the latent factor matrices $\matn{A}\in\Real^{I\times R}$ (i.e. the ground truth) were drawn from i.i.d. standard normal distributions, $n=1,2,\ldots,N$. Then independent Gaussian noise (with SNR=10dB) was added to $\tensor{Y}^*$ to yield noisy observation tensor \tensor{Y}.  We tested three different combinations of $N,\; I,$ and $R$, as shown in TABLE \ref{tab:simuCP}. The algorithms CP-ALS, PARAFAC, which are the standard algorithms included in the Tensor Toolbox 2.5 \cite{KoldaTensorToolbox} and the $N$-way toolbox \cite{nwaytoolbox}, respectively, were compared as the baseline. The FastCPALS algorithm  \cite{PhanCPgrad} was also compared. For all the algorithm, the maximum iteration number was set as 1000, and a stopping criterion $|\text{Fit}(\tensor{Y},\tensor{\hat{Y}}_t)-\text{Fit}(\tensor{Y},\tensor{\hat{Y}}_{t+1})|<10^{-6}$ was applied, where  $\tensor{\hat{Y}}_t$ is an estimation of $\tensor{Y}^*$ obtained after the $t$th iteration. The mean values and standard variations of Fit($\tensor{Y}^*,\tensor{\hat{Y}}$) over 20 Monte-Carlo runs are detailed in TABLE \ref{tab:simuCP}, along with the averaged time consumed by each algorithm. From the table, we see that the FFCP algorithm is significantly faster than the other algorithms, but without  significant deterioration of performance, although it works on a compressed tensor instead of the original tensor.  We also tested the PARACOMP \cite{paracomp} method which was specially proposed for big tensor decompositions. It turned out that the PARACOMP was sensitive to noise: in noise-free cases, it often achieved perfect Fit of 1; otherwise when the noise was 10dB, it generally obtained the Fit lower than 0.3.

In \figurename \ref{fig:scalability} we demonstrated the scalability of the proposed FFCP method where the RandTucker2i was used to compress the data with $p=15$. From \figurename \ref{fig:scalability}(a) we can see  that the FFCP shows very promising scalability: its time consumption  increases almost linearly with respect to the scale of the input data on condition that the rank is fixed; on the contrast, the time consumed by the other algorithms increased nearly exponentially. Note that both the FFCP and the PARACOMP method \cite{paracomp} were designed mainly for big data analytics, and this comparison did not precisely reflect their actual time consumption difference in a distributed computational environment. However, it did show that the proposed method was more efficient than the PARACOMP, which can also be inferred from the obvious facts: 1) the PARACOMP compresses the observation big tensor at least $P$ times ($P\ge {I}/{R}$) whereas in our RandTucker2i based FFCP we only need to do it twice (or even once when RandTucker is applied); 2) although both of them need to perform CPD on small compressed tensors, however, the PARACOMP method needs to do this $P$ times whereas the FFCP needs to do it only once as the final step. Particularly, in the PARACOMP any failure of the CPDs of $P$ sub-tensors will degenerate the final results. \figurename \ref{fig:scalability}(b) shows the scalability of the FFCP provided that the given big tensor was in its Tucker format. To such scale of tensor, if we use the original full tensor directly, all the algorithms ran out of memory. But if we use its Tucker compression, the FFCP can perform  constrained/unconstrained CPD efficiently even on our desktop computer.

{\bf Simulations on Nonnegative CPD.}  In this experiments the data were generated as the above, except that the latent factor matrices were generated from i.i.d. exponential distributions with the mean parameter $\mu=10$. Then we randomly selected 10\% entries of each factor to be zeros to generate sparse factor matrices. As the PARAFAC algorithm was very slow, we only compared the CP-NMU algorithm, which is a standard nonnegative CPD algorithm based on the multiplicative update rule included in the Tensor Toolbox \cite{KoldaTensorToolbox}, and the CP-HALS algorithm \cite{TensorHALS2009} that is based on the Hierarchical  ALS (HALS) method. For fair comparison, their corresponding implementations of FFCP were compared. The results averaged over 20 Monte-Carlo runs were detailed in TABLE \ref{tab:simuNTF}. From the TABLE, FFCP implementations were not only significantly faster than their competitors, but also often achieved higher Fit values. To the best of our knowledge, the FFCP is one of the most efficient algorithms for nonnegative CPD.

To demonstrate the high efficiency of the FFCP algorithm, in this simulation the entries of $\matn{A}\in\Real^{10,000\times 10}$, $n=1,2,3,4$, were drawn from i.i.d. exponential distributions with the mean parameter $\mu=10$. Again, we randomly selected 20\% entries of \matn{A} to be zeros to obtain sparse factor matrices. By using \eqref{eq:Polyadic} and \matn{A} generated above, a big 4th-order tensor $\tensor{Y}\in\Real^{10,000\times 10,000\times10,000\times10,000}$ with $R=10$ was obtained, which contains a total of $10^{16}$ entries in its full format. It is certainly impossible to decompose such a huge tensor in a standard desktop computer. Suppose we have already obtained its Tucker approximation $\tensor{Y}_{\text{T}}=\tenfactors{\tensor{G};\matn[1]{A},\matn[2]{A},\matn[3]{A},\matn[4]{A}}\approx\tensor{Y}$. Due to the rotation ambiguity of Tucker decompositions, we know theoretically $\matn{\tilde{A}}=\matn{A}\matn{Q}$, where \matn{Q} was not necessarily to be nonnegative. Generally applying CPD on the core tensor \tensor{G} directly cannot guarantee the nonnegativity of factor matrices. Instead, we applied the FFCP algorithm on the Tucker compressed tensor $\tensor{Y}_{\text{T}}$ to extract nonnegative factor matrices. In a typical run the FFCP-HALS consumed only 77.4 seconds, and the recovered components were with the signal-to-interference ratio\footnote{The performance index SIR is defined as $\text{SIR}(\matn{a}_r,\matn{\hat{a}}_r)=20\log{\frob[2]{\matn{a}_r-\matn{\hat{a}}_r}}/{\frob[2]{\matn{a}_r}}$, where $\matn{a}_r$ is an estimate of $\matn{a}_r$, both $\matn{a}_r$ and $\matn{a}_r$ are normalized to be zero-mean and unit-variance.} (SIR) higher than 200dB, which almost perfectly match the true components. Although this example only shows an ideal noise-free case, it illustrates how the FFCP algorithm can be applied to perform nonnegative CPD of huge tensors, as long as their good low-rank approximations can be obtained.


\subsection{The COIL-100 Images Clustering Analysis}
\begin{figure}[!t]
\centerline{
    \includegraphics[width=\linewidth]{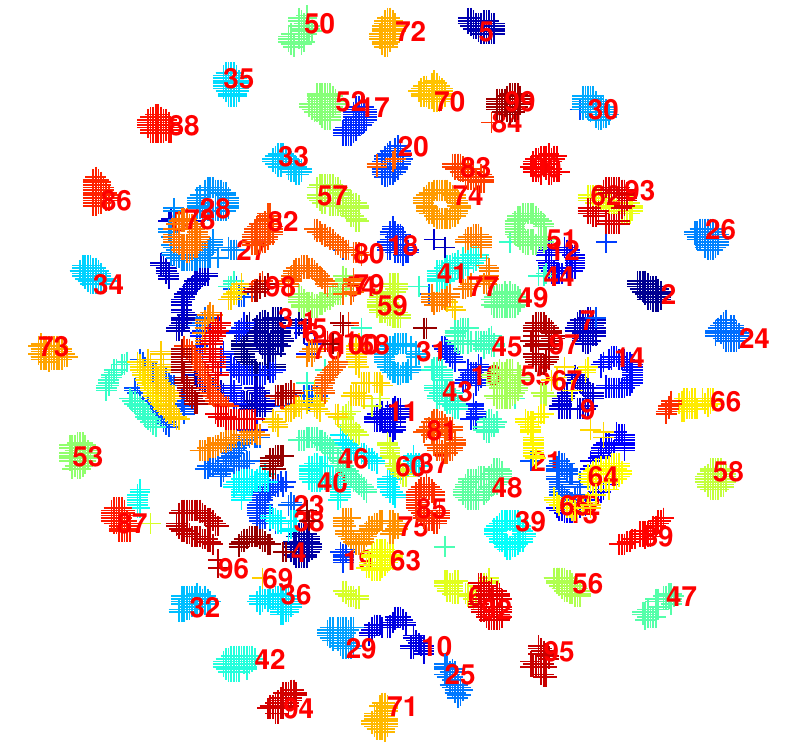}
}
\caption{Visualization of the clustering results obtained by applying distributed RandTucker2i on the COIL100 image database. The two tSNE components of \matn[4]{A} were used for visualization.}
\label{fig:COIL100clusters}
\end{figure}

In this experiment, we applied the proposed NTD algorithms to the clustering analysis of the objects selected from the Columbia Object Image Library (COIL-100).  The COIL-100 database consists of 7,200 images of 100 objects, each of which has 72 images taken from different poses, hence these images naturally form a $128\times 128 \times 3 \times 7200$ tensor \tensor{Y}. The data tensor was decomposed by using tensor decomposition algorithms and then used the two t-Distributed Stochastic Neighbor Embedding (t-SNE) components of \matn[4]{A} as features to cluster and visualize the results. The K-means algorithm was used to cluster the objects and the accuracy of clustering was defined in  \cite{GNMF2011PAMI}.To justify the proposed distributed RandTucker algorithm, we used the MATLAB parallel computing toolbox which allows us to use the full processing power of multicore desktops by executing applications on workers (MATLAB computational engines) that run locally. The tensor \tensor{Y} was partitioned into 4 sub-tensors with the size of $128\times128\times3\times1800$ and was distributed on 4 workers. As the baseline the Tucker-ALS algorithm included in \cite{KoldaTensorToolbox} was compared, and the multilinear rank was empirically set as $\mat{R}=(20,\,20,\,3,\,20)$. For the RandTucker method $\mat{\tilde{R}}=(30,\,30,\,3,\,50)$ was set. Both of these two algorithms used 2 iterations. Their performances were detailed in TABLE \ref{tab:simuCOIL}. It can be seen that no significant deterioration was observed between the randomized methods and the deterministic methods, which suggested that the distributed RandTucker algorithm can be applied to really large scale problems. The clustering results obtained by the RandTucker2i \figurename are visualized in \figurename \ref{fig:COIL100clusters} by using the two tSNE components of \matn[4]{A}.

\begin{table}
  \caption{Comparison between the algorithms when they were applied to clustering analysis of COIL-100 objects.}
  \label{tab:simuCOIL}
  \centering
  \begin{tabular}{ r  c c } \hline \hline
  Algorithm & Fit (\%) & Accuracy (\%) \\ \hline
 Tucker-ALS & 70.3 & 69.2 \\
 RandTucker2i & 69.1 & 68.4 \\
  RandTucker & 68.6 & 69.8 \\
 FFCP(Tucker-ALS) & 67.1 & 71.3 \\
 FFCP (RandTucker2i) & 64.7 & 68.0 \\ 
 \hline \hline
  \end{tabular}
  \end{table}


\subsection{UCF Sports Action Videos Analysis}
This data set consists of a set of videos collected from various sports which are typically featured on broadcast television channels \cite{UCFSports} (see \url{http://crcv.ucf.edu/data/UCF_Sports_Action.php}). To avoid rescaling the videos and keep all the videos in the same specification we selected 5 actions of Diving, Riding-Horse, Run, Swing, and Walk-Front, each of which contains 10 video clips with the size of $404\times720\times3\times35$ (35 is the number of frames). We concatenated these videos to form a 5th-order tensor. This tensor has about $1.6\times 10^9$ entries and more than 99\% of them are nonzero. When it is saved in double format it consumes about 11.7 GB space. We used tensor decomposition methods to decompose this tensor to extract low-dimensional features for further clustering analysis. When we applied the CP-ALS method included in the Tensor Toolbox to decompose 30 videos, it consumed about 4 hours for only one iteration; when 40 videos were considered, the CP-ALS ran out of memory. We configured a small cluster with 3 computers (a desktop computer plus two laptop computers (each of which has a i7 Intel dual core CPU, 8 GB physical memory running Windows 7 and MATLAB 2012b)). In the first stage the RandTucker2i was applied with $\mat{\tilde{R}}=[35,35,3,10,20]$. It consumed about 104 seconds to obtain the Tucker approximation of the original tensor, with the Fit of 72.3\% to the observation data. Then the FFCP was applied using a desktop computer with $R=10$. In this case only 7.5 seconds were required without nonnegativity constraints; and 7 seconds were consumed when nonnegativity constraints were imposed. Finally we used \matn[5]{A} as features to cluster the videos and achieved the accuracy of 58\% in a typical run. This experiments showed that the proposed algorithms can be applied to decompose very large scale tensors by using a distributed computing cluster.

 \begin{figure}[!t]
\centerline{
    \includegraphics[width=\linewidth]{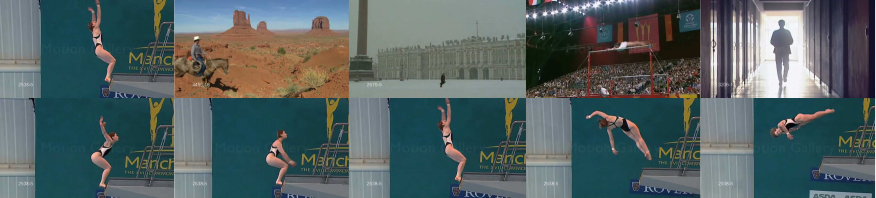}
}
\caption{The UCF Sports Action Videos. The first row shows a snapshot of each action, and the second row shows 5 snapshots from a video of diving.}
\label{fig:UCFSports}
\end{figure}

\section{Conclusion}
\label{sec:conclusion}
Tensor decompositions are particularly promising for big multiway data analytics, as it is able to discover more complex internal structures and correlations of data gathered from multiple aspects, devices, trials, and subjects, etc. While most existing tensor decomposition methods are not designed to meet the major challenges posed by big data, in this paper we proposed a unified framework for the decomposition of big tensors with relatively low multilinear rank. The key idea is that based on its Tucker approximation of a big tensor, a large set of unconstrained/constrained tensor decompositions can be efficiently performed. We developed a flexible fast algorithm for the CP decomposition of tensors using their Tucker representation, and a distributed randomized Tucker decomposition approach which is particularly suitable for big tensor analysis. Simulation results confirmed that the proposed algorithms are highly efficient and scalable. 

While the randomized Tucker decomposition approach showed its close-to-optimal performance, the corresponding error bounds need to be investigated theoretically in the future. We also believe many randomized (or related) approaches proposed for low-rank approximation of matrices \cite{siam_probLowRank,LRApIEEE} could be applied for Tucker decompositions, however, a comparative study in the tensor decomposition scenarios is desired. Particularly, if we have efficient and robust unconstrained CP decomposition algorithms, we can use the CP model to compress the raw data instead of the Tucker model which suffers from the curse of dimensionality for very high-order tensors. Hence developing more reliable and robust CP decomposition approaches is another important task. All these open issues deserve further study in the near future.

\bibliographystyle{IEEEtran}



\end{document}